\documentclass[11pt]{article}
\usepackage{amsmath, amsfonts, amsthm, amssymb, verbatim}
\usepackage{graphicx,color}
\usepackage{amstext}
\usepackage{amssymb,mathrsfs}

\numberwithin{equation}{section}
\newcommand{\zo}{{\mathbb Z ^2_0}}
\newcommand{\zp}{{\mathbb Z ^2_+}}

\newtheorem{theorem}{Theorem}[section]
\newtheorem{proposition}[theorem]{Proposition}
 \newtheorem{remark}[theorem]{Remark}
\newtheorem{lemma}[theorem]{Lemma}
\newtheorem{corollary}[theorem]{Corollary}

\begin{document}

\title{2D Navier-Stokes equation with cylindrical fractional Brownian noise}

\author{Benedetta Ferrario\thanks{Dipartimento di Matematica ''F. Casorati'', Universit\`{a} di Pavia,
Italy, E-mail: \textsl{benedetta.ferrario@unipv.it}}, Christian
Olivera\thanks{Departamento de Matem\'{a}tica, Universidade Estadual de
Campinas, Brazil. E-mail: \textsl{colivera@ime.unicamp.br}. }}

\date{\today}

\maketitle

%%%%%%%%%%%%%%%%%%%%%%%%%%%%%%%%%%%%%%%%%%%%%%%%%%%%%%%%%%%%%%%%%%%%%%%%%%%%
\begin{abstract}
We consider the Navier-Stokes equation on the 2D torus, with a stochastic forcing term 
which is a cylindrical fractional Wiener noise of Hurst parameter $H$. Following 
\cite{DPD,AF} which dealt with the case $H=\frac12$, 
we prove a local existence and uniqueness result when 
$\frac 7{16}< H<\frac 12$ and a global existence and uniqueness result
 when $\frac 12<H<1$.
\end{abstract}
%%%%%%%%%%%%%%%%%%%%%%%%%%%%%%%%%%%%%%%%%%%%%%%%%%%%%%%%%%%%%%%%%%%%%%%%%%%%


\vspace{0.3cm} \noindent {\bf MSC2010 subject classification:} 60H15, %SPDE
 35R60, %Partial differential equations with randomness, stochastic partial differential equations 
 60H30, %Applications of stochastic analysis (to PDE, etc.)
76D05. % 	Navier-Stokes equations

%%%%%%%%%%%%%%%%%%%%%%%%%%%%%%%%%%%%%%%%%%%%%%%%%%%%
\section {Introduction} \label{Intro}
An incompressible fluid flow dynamics is described by the so-called incompressible Navier-Stokes
equations. 
In this paper,  we consider the Navier-Stokes equations on the torus, 
i.e. we work on the square 
$\mathbb T=[0,2\pi]^2$ with periodic boundary conditions;  we add a stochastic forcing term.
These are the equations
\begin{equation}\label{eq-initial}
 \left \{
\begin{aligned}
    &\partial_t v =  \nu \Delta v - (v\cdot \nabla) v  - \nabla p + \partial_{t} W^{H}
    \\[2pt]
    &\text{div }v=0 \\[2pt]
    &v|_{t=0}=  v_{0} 
\end{aligned}
\right .
\end{equation}
where for $t\in [0,T]$ and $\xi \in \mathbb T$, $v=v(t,\xi)$ is the vector velocity, $p=p(t,\xi)$ the scalar pressure, $\nu > 0$ the viscosity coefficient
 and $W^H=W^H(t,\xi)$ a cylindrical fractional Brownian process. 

Stochastic perturbations in the equations of motions are commonly used to model small perturbations
(numerical, empirical, and physical uncertainties) or thermodynamic fluctuations
present in fluid flows. We refer to the lecture notes by Flandoli \cite{Fla}, the monograph
of Kuksin and Shirikyan \cite{Kuk} as well as the references cited therein for a recent overview.

Different noise terms have been considered so far. 
The contribution of this paper is to study equation \eqref{NS} with a cylindrical fraction Brownian motion $W^H$ for $H \neq \frac 12$. 
Indeed the case $H=\frac 12$ has been studied in \cite{AC,AF,DPD,Deb}.
Let us point out that with a coloured (not cylindrical) noise, the analysis of equations \eqref{eq-initial} is easier; some results on a more general bidimensional domain can be found in \cite{FSV}. Moreover, 
the bigger is $H$ the more regular is the fractional Brownian motion. Hence it is worth to ask if the cylindrical fractional Brownian motion
with $H<\frac 12$ can be considered; in addition also the analysis for $H>\frac 12$ is interesting in order  to compare the results for different values of the Hurst parameter. 
In this paper 
we shall prove local existence and uniqueness of solutions for $\frac 7{16}<H<\frac 12$ 
and  global existence and uniqueness for $\frac 12<H<1$. 

As far as the contents of the paper are concerned, in Section 2 we introduce the mathematical setting, 
in Section \ref{sectionStokes} we analyze the linear Stokes problem, whereas Section \ref{sectionB}
analyzes the bilinear term  and Section \ref{sectionNS} the  Navier-Stokes problem.
In the Appendices we present some proofs.

\section{Mathematical setting}
In this section we introduce the basic tools.

\subsection{The spaces}
For a complex number $b=\Re b+i\Im b$ we denote by $\overline b$ the complex conjugate 
($\overline b=\Re B-i \Im b$)
and by 
$|b|$ the absolute value ($|b|=\sqrt{(\Re b)^2+(\Im b)^2}$). 

We consider subspaces of $\mathbb Z^2$:
\[
\mathbb Z^2_0
=\{k=(k^{(1)},k^{(2)})\in \mathbb Z^2: k\neq 0\}
\]
\[
\zp
=\{k=(k^{(1)},k^{(2)})\in \mathbb Z^2_0: k^{(1)}>0\}\cup \{k=(k^{(1)},k^{(2)})\in \mathbb Z^2_0: k^{(1)}=0, k^{(2)}>0\}
\]
and 
\[
\mathbb Z^2_-=\mathbb Z^2_0\setminus \mathbb Z^2_+
\]
When $k=(k^{(1)},k^{(2)})\in \mathbb Z^2$, we denote by $|k|$ 
the absolute value ($|k|=\sqrt {(k^{(1)})^2+(k^{(2)})^2}$).

We consider the separable Hilbert space $\mathcal H^0$  which is the $L^2$-closure of the space of smooth  vectors which are
 periodic, zero mean value 
and divergence free.
Let  $\{h_k\}_{k}$ be the basis for $\mathcal H^0$, given by
$h_k(\xi)=\frac 1{2\pi} \frac {k^\perp}{|k|}e^{i k \cdot \xi}$ for $k \in \mathbb Z^2_0$ 
and $\xi \in \mathbb T$. 
Notice that, for any $k \in \mathbb Z^2_0$,
$h_{-k}(\xi)=-\overline{h}_k(\xi)$  
 and $\Delta h_k=-|k|^2 h_k$. Therefore
\[
\mathcal H^0=\{v(\xi)=\sum_{k\in \zo} v_k h_k(\xi): v_{-k}=-\overline v_k\ \forall k, 
\sum_{k\in\zo} |v_k|^2<\infty\}
\]
Notice that the complex coefficients $v_k$ must satisfy 
$v_{-k}=-\overline v_k$ in order to get a real vector $v$.

More generally, for $r\in \mathbb R$  we define 
\[
\mathcal H^r=\{v(\xi)=\sum_{k\in \zo} v_k h_k(\xi): v_{-k}=-\overline v_k\ \forall k, 
\sum_{k\in\zo} |k|^{2r}|v_k|^2<\infty\}.
\]
This is a Hilbert space with scalar product
\[
(u,v)_{\mathcal H^r}=\sum_{k\in \zo} |k|^{2r}u_k \overline v_k .
\]

Following \cite{BL},
we define the periodic divergence-free vector Sobolev spaces ($r\in
\mathbb R, 1\le p\le \infty$)
$$
 {\mathcal H}^r_p=\{ \;
             v =\sum_{k \in \zo} v_k h_k:
                       \sum_{k\in\zo} v_k |k|^r h_k \in [L^p(\mathbb T)]^2\ \}
$$
and the  periodic divergence-free  vector Besov spaces as  real interpolation
spaces
\begin{equation*}
\begin{split}
 \mathcal B^{r}_{p\,q}=(\mathcal H^{r_0}_p,
                          \mathcal H^{r_1}_p)_{\theta,q}
   , \qquad &r \in \mathbb R, 1\leq p,q\le \infty
\\
 & r=(1-\theta)r_0+\theta r_1,
     \qquad
 0<\theta<1
\\
\end{split}
\end{equation*}
In particular $\mathcal B^r_{2\,2}=\mathcal H^r_2 =\mathcal H^r$. Moreover (see \cite{BL})
\[
\begin{split}
\|v\|_{\mathcal B^s_{p\, q_1}} &\le \|v\|_{\mathcal B^s_{p\, q_2}} \qquad \text{ for } q_2\le q_1\\
\|v\|_{\mathcal B^{s_1}_{p\, q}} &\le \|v\|_{\mathcal B^{s_2}_{p\, q}} \qquad \text{ for } s_1\le s_2\\
\|v\|_{\mathcal B^{s_1}_{p_1\, q}} &\le C\|v\|_{\mathcal B^{s_2}_{p_2\, q}} \qquad \text{ for } s_1-\frac 2{p_1}=s_2-\frac2{p_2} 
\end{split}
\]
Here $C$ is a generic constant. We make the convention to denote
  different constants by the same symbol $C$, unless we want to mark them
  for further reference.

One interesting result in Besov spaces is given by the following estimate 
 of Chemin (see Corollary 1.3.1 in \cite{Ch}):
\begin{equation}\label{chemin}
\|v_1 v_2\|_{\mathcal B^s_{pq}}\le 
\frac{C^{s_1+s_2}}{s_1+s_2}\|v_1\|_{\mathcal B^{s_1}_{pq}} \|v_2\|_{\mathcal B^{s_2}_{pq}}
\end{equation}
if 
\[
s_1+s_2>0, \qquad 
s_1<\frac 2p, \qquad s_2<\frac 2p, \qquad s=s_1+s_2-\frac 2p 
\]
and $p,q \in [1,\infty]$.

%%%%%%%%%%%%%%%%%%%%%%%%%%%%%
\subsection{The abstract equation}
Let us consider a unitary viscosity $\nu=1$ in system \eqref{eq-initial}. Then we
 write the evolution  in abstract form as
\begin{equation}\label{NS}
dv(t)=Av(t)\ dt -B(v(t),v(t))\ dt +dw^H(t)
\end{equation}
with  the operators formally defined as $A=\Delta$ and $B(u,v)=P[(u\cdot\nabla)v]$, where 
$P$ is the projector operator onto the space of divergence free vector fields.
We can represent the stochastic forcing term as 
\begin{equation}\label{noise}
w^H(t,\xi)=\sum_{k\in\zo}h_k(\xi) b^H_k(t), \qquad  (t, \xi) \in \mathbb R \times\mathbb T
\end{equation}
where $\{b^H_k\}_{k \in \mathbb Z^2_+}$ is a sequence of i.i.d. complex fractional Brownian processes
defined on a complete probability space $(\Omega,\mathbb F, \mathbb P)$ with filtration 
$\{\mathbb F_t\}_{\{t \in \mathbb R\}}$ and $b^H_{-k}=\overline{b^H_k}$ for all $k \in \zp$. 
We denote by $\mathbb E$ the mathematical expectation with respect to $\mathbb P$. 
This means that
$b^H_k(t)=\Re b_k^H(t)+i\Im  b_k^H(t)$ and $\{\Re b^H_k, \Im b_k^H(t)\}_{k \in \mathbb Z^2_+}$ 
is a sequence of i.i.d. standard real fractional Brownian processes (fBm) with Hurst parameter $H$.
Each element of the sequence is a centered Gaussian process whose covariance 
is
\[
C(t,s)=\frac 12 (|t|^{2H}+|s|^{2H}-|t-s|^{2H})
\]
$w^H$ is called an $\mathcal H^0$-cylindrical fractional Brownian motion and one can prove that the series \eqref{noise} converges in any space $U$  with 
continuous embedding $\mathcal H^0\subset U$ of Hilbert-Schmidt type.

Now let us define rigorously the operators $A$ and $B$.

The Stokes operator $A$, as a linear operator in $\mathcal B^{s}_{pq}$ with domain $\mathcal B^{s+2}_{pq}$,
 generates an analytic semigroup $\{e^{tA}\}_{t\ge 0}$ in $\mathcal B^{s}_{pq}$ and
\begin{equation}\label{stima-expA}
\|e^{tA}v\|_{\mathcal B^{s_1}_{pq}} \le \frac C{t^{\frac{s_1-s_2}2}} \|v\|_{\mathcal B^{s_2}_{pq}}
\end{equation}
for any $t\ge 0$, $s_1>s_2$.

As far as the bilinear term $B(u,v)=P[(u\cdot \nabla)v]$ is considered, we recall some basic properties
(see \cite{Temam}). Let $\langle\cdot,\cdot\rangle$ denote the 
$\mathcal H^{-r}-\mathcal H^r$ duality bracket. One checks by integrations by parts that 
\begin{equation}\label{tril}
\langle B(u_1,u_2),u_3\rangle=-\langle B(u_1,u_3),u_2\rangle
\end{equation}
and taking $u_2=u_3$
\begin{equation}\label{Bnull}
\langle B(u_1,u_2),u_2\rangle=0 .
\end{equation}
These relationships are true with regular entries and  then are extended to 
more general vectors by density. 

A basic estimate is (see \cite{Giga}, Lemma 2.2)
\begin{equation}\label{stimeGiga}
\|B(u,v)\|_{\mathcal H^{-\delta}}\le C\|u\|_{\mathcal H^{\theta}} \|v\|_{\mathcal H^{\rho}}
\end{equation}
when 
\[
0\le \delta<2,\qquad \rho>0,\qquad \theta>0
\]
\[ 
\rho+\delta>1,\qquad
\theta+\rho+\delta\ge 2
\]
Other estimates have been given before in \eqref{chemin}; indeed, by the divergence free condition we have
$B(u,v)=P[\text{div}\ (u \otimes v)]$; therefore $\|B(u,v)\|_{\mathcal B^s_{pq}}\le \|u \otimes v\|_{\mathcal B^{s+1}_{pq}}$.

Moreover, as done in  \cite{AF}, we can develop the bilinear term in Fourier series.
Given $v=\sum_{l\in \zo} v_l h_l$ and $u=\sum_{h \in \zo} u_h h_h$, 
we have formally
\[
\begin{split}
B(u,v)&=iP\sum_{h,l\in \mathbb Z^2_0}u_h \frac {h^\perp\cdot l}{|h|}\frac {e^{i h \cdot \xi}}{2\pi}
   v_l \frac {l^\perp}{|l|}\frac {e^{i l \cdot \xi}} {2\pi}
\\&=
iP\sum_{k\in \zo}\left(\sum_{\substack{h \in \zo\\ h\neq k}}
         \frac{h^\perp\cdot k}{2\pi |h||k-h|}u_hv_{k-h} (k-h)^\perp \right) 
\frac {e^{i k \cdot \xi}} {2\pi}
\end{split}
\]
Using that the projector $P$ acts on the $k$-th component as 
$P_k a=\frac{a\cdot k^\perp}{|k|^2}k^\perp$, we get
\[
B(u,v)
= i\sum_{k\in \zo}\left(\sum_{\substack{h \in \zo\\ h\neq k}}
         \frac{h^\perp\cdot k}{2\pi |h||k-h|}u_hv_{k-h} \frac{(k-h)^\perp \cdot k^\perp}{|k|} \right) 
k^\perp \frac {e^{i k \cdot \xi}} {2\pi|k|}
\]

Summing up, the bilinear term can be written in Fourier series as
\begin{equation}\label{B-inFourier}
B(u,v)=\sum_{k \in \zo} B_k(u,v)h_k
\end{equation}
with 
\begin{equation}\label{componentiBk}
\begin{split}
  B_k(u,v)&=i  \sum_{\substack{h \in \zo\\ h \neq k}} \gamma_{h,k} u_h v_{k-h}
\\
  \gamma_{h,k}\;&=\frac{1}{2\pi}\frac{(h^\perp \cdot k)([k-h]\cdot k)}{|h||k-h||k|}
\\
\end{split}
\end{equation}
Notice that $\overline B_k=-B_{-k}$.
The convergence of the series \eqref{B-inFourier} will be analysed in the next section.

Our aim is to study equation \eqref{NS} for $H \neq \frac 12$. 
Indeed the case $H=\frac 12$ has been studied in \cite{AC,AF,DPD,Deb}:
Da Prato and Debussche proved the existence of a strong mild solution for 
$\mu$-a.e. initial condition (where $\mu$ is the Gibbs measure of the enstrophy, introduced in 
\cite{AFHK} which is an invariant measure for equation \eqref{NS}), 
whereas Albeverio and Ferrario proved pathwise uniqueness of these solutions.

We shall prove  a local existence and uniqueness result for  $\frac 7{16}<H<\frac 12$ 
  and a global existence  and uniqueness result
for   $H>\frac 12$. This latter result  improves that  of
\cite{FSV}; indeed, the case of cylindrical fBm is 
included in \cite{FSV} but only for $H>\frac 34$ (see Theorem 5.1 and Corollary 4.3 there).
By the way, there are other  differences with respect to \cite{FSV}: 
in \cite{FSV} the spatial domain is not the torus but a generic smooth bounded
subset $D$ of $\mathbb R^2$ (and the Dirichlet boundary condition is assumed) and the solution is a 
process with values in $L^4$ in time and space, whereas our solution is more regular
since a.a. paths are at least in
$C([0,T];\mathcal H^{\frac 12})$ (see next Theorem \ref{th-global-ex} with $H>\frac 34$)
and one knows that  $\mathcal H^{\frac 12}\subset L^4(D)$.

\medskip
In order to analyze equation \eqref{NS} we introduce as in 
\cite{DPD} two subproblems: the linear Stokes equation
\[
dz(t)=Az(t)\ dt +dw^H(t)
\]
and the equation for $u=v-z$
\[
\frac{du}{dt}(t)=Au(t)-B(u(t),u(t))-B(u(t),z(t))-B(z(t),u(t))-B(z(t),z(t))
\]
which is a Navier-Stokes type equation with random coefficients.

First we deal with the linear problem for $z$, then we define the bilinear 
term $B(z,z)$ a.s.  as in \cite{DPD,AF} 
and finally face the nonlinear equation for $u$. At the end we recover the 
existence result for $v$ from the representation $v=z+u$.

\section{The  Stokes equation}\label{sectionStokes}
If we neglect the bilinear term in \eqref{NS}, we obtain the linear Stokes equation
\begin{equation}\label{eq-lin}
dv(t)=Av(t)\ dt+dw^H(t).
\end{equation}
We consider its stationary mild solution; this is the process
\begin{equation}\label{OU}
z(t)=\int_{-\infty}^t e^{(t-s)A}dw^H(s)
\end{equation}
We can write
\begin{equation}\label{z-in-Fourier}
\begin{split}
z(t)(\xi)&=\sum_{k\in \zo} h_k(\xi) \int_{-\infty}^t e^{-|k|^2(t-s)}db_k^H(s)
\\&
\equiv
2\sum_{k\in\zp}\frac {k^\perp}{2\pi|k|}\cos(k\cdot\xi) \int_{-\infty}^t e^{-|k|^2(t-s)}d\Re b_k^H(s)
\\&\qquad -
2\sum_{k\in\zp}\frac {k^\perp}{2\pi|k|}\sin(k\cdot\xi) \int_{-\infty}^t e^{-|k|^2(t-s)}d\Im b_k^H(s)
\end{split}
\end{equation}
First, we provide a result for  each stochastic convolution integral 
appearing in the Fourier series representation.
%LEMMA
\begin{lemma}
Let $\lambda>0$ and $b^H$ be a real fBm of Hurst parameter $H\in (0,1)$.
Then
\[
\int_{-\infty}^t e^{-\lambda(t-s)}db^H(s), \qquad t \in \mathbb R
\]
is a stationary centered Gaussian process whose variance is
\[
C_H\lambda^{-2H}
\]
where $C_H$ is the positive constant given in \eqref{costanteCH}.
\end{lemma}
\proof
 Following the proof of Lemma 4.1 in \cite{FSV} we have that the random variables
\[
\int_{-\infty}^t e^{-\lambda (t-s)}db^H(s)
\]
and
\[
\int_0^{+\infty} e^{-\lambda r}db^H(r)
\]
have the same law.
Moreover, by self-similarity of the fBm,
the latter random variable has the same law as 
\[
\lambda^{-H} \int_0^{+\infty} e^{-r}db^H(r).
\]

Therefore 
\[
\mathbb E\left(\int_{-\infty}^t e^{-\lambda (t-s)}db^H(s)\right)^2
= \lambda^{-2H}
\mathbb E\left( \int_0^{+\infty} e^{-r}db^H(r)\right)^2
\]
We estimate $\mathbb E\left( \int_0^{+\infty} e^{-r}db^H(r)\right)^2$
using the representation
\[
\int_0^{+\infty} e^{-r}db^H(r)=\int_0^{+\infty} e^{-r}b^H(r) dr.
\]
This comes from the formula on a finite time interval
\[
\int_0^{T} e^{-r}db^H(r)=e^{-T}b^H(T)+\int_0^{T} e^{-r}b^H(r) dr
\]
and the fact that by the law of iterated logarithm (see \cite{BHOZ}) we get 
\[
\lim_{T\to +\infty}|e^{-T}b^H(T)|=0 \qquad \mathbb P-a.s.
\]

Hence
\[\begin{split}
\mathbb E\left( \int_0^{+\infty} e^{-r}db^H(r)\right)^2
&=
\mathbb E\left( \int_0^{+\infty} e^{-r}b^H(r)dr\right)^2
\\&=
\int_0^{+\infty}\int_0^{+\infty}e^{-r}e^{-s} \frac {r^{2H}+s^{2H}-|r-s|^{2H}}2 \ dr \ ds
\end{split}\]
By elementary calculations one shows that the latter integral is finite.
We set
\begin{equation}\label{costanteCH}
C_H=\int_0^{+\infty}\int_0^{+\infty}e^{-r}e^{-s} \frac {r^{2H}+s^{2H}-|r-s|^{2H}}2 \ dr \ ds.
\end{equation}
\hfill\qed%\smartqed

Now we come back to the stationary process $z$ given in \eqref{OU}. We have
the following result
\begin{proposition}\label{pro-zLm}
For any $r<2(H-\frac 12)$   we have
\[
z\in C(\mathbb R;\mathcal H^r) \qquad \mathbb P-a.s.
\]
\end{proposition}
\proof
First we show that for any fixed time, the random variable $z(t)\in \mathcal H^r$, $\mathbb P$-a.s.
Indeed, using \eqref{z-in-Fourier} and the previous Lemma we have
\[\begin{split}
\mathbb E\big[\|z(t)\|_{\mathcal H^r}^2\big]
&=\sum_{k\in\zo} |k|^{2r}\left|\int_{-\infty}^te^{-|k|^2(t-s)}db_k^H(s)\right|^2
\\&
=\sum_{k\in\zo} |k|^{2r}2\frac{C_H}{|k|^{4H}}
\end{split}\]
The latter series is convergent for $4H-2r>2$, i.e. $r<2(H-\frac 12)$.

It follows that for any finite $m\ge 1$ we have
$z\in L^m_{\text{loc}}(\mathbb R;\mathcal H^r)$, $\mathbb P$-a.s.. Indeed,
$z(t)$ is a Gaussian random variable; so  all the moments are finite, i.e. 
for any $m\ge 2$ there exists a finite constant $e_m$ such that
\[
\mathbb E[\|z(t)\|_{\mathcal H^r}^m]=e_m
\]
for any $t$.
Moreover,  the process $z$ is a stationary process and by interchanging the integrals, 
for any $T_1<T_2$ we get
\[
\mathbb E\left[\int_{T_0}^{T_1} \|z(t)\|_{\mathcal H^r}^m dt \right]=
\int_{T_0}^{T_1}\mathbb E \left[\|z(t)\|_{\mathcal H^r}^m \right] dt=e_m(T_1-T_0)<\infty
\]
Since the expectation is finite, then $\int_{T_0}^{T_1} \|z(t)\|_{\mathcal H^r}^m dt<\infty$, $\mathbb P$-a.s.

The continuity in time of the trajectories
has been proved  in \cite{DDM2002} when $H>\frac 12$ and in \cite{DDM2006}   when $H<\frac 12$.
\hfill\qed%\smartqed

\begin{remark}
We see that when $H\le\frac 12$, the process $z$ at any fixed time takes values in a distributional space. 
This is the source of the difficulty in our problem.
\end{remark}

\begin{remark}
From the proof of Proposition \ref{pro-zLm}, we obtain that the process $z$ is a stationary process 
and for any time $t$ the law of $z(t)$ is the  centered Gaussian measure
$\mu^H\sim\mathcal N(0,C_H(-A)^{-2H})$. 
More precisely,  we assign the measure $\mu^H$ on the  sequences $\{(\Re v_k,\Im v_k)\}_{k \in \zp}$
as
\begin{equation}\label{measure}
\mu^H=\otimes_{k \in \zp}\mu^H_k
\end{equation}
with 
\[
d\mu^H_k(x,y)=\frac{|k|^{4H}}{2\pi C_H}e^{-\frac{|k|^{4H}}{2C_H}(x^2+y^2)}\ dx \ dy
\]

When we identify the space $\mathcal H^r$ with that of the sequences 
$\{(\Re v_k,\Im v_k)\}_{k \in \zp}$ such that $\sum_{k\in\zp} |k|^{2r}[(\Re v_k)^2+(\Im v_k)^2]<\infty$,
we get
$\mu^H(\mathcal H^r)=1$ for any $r<2(H-\frac 12)$ and  
$\mu^H(\mathcal H^r)=0$ for any $r\ge2(H-\frac 12)$
(see \cite{Kuo}). Similarly, $\mu^H(\mathcal B^r_{pq})=1$ for any $r<2(H-\frac 12)$ and  
$\mu^H(\mathcal B^r_{pq})=0$ for any $r\ge2(H-\frac 12)$.
\end{remark}

We finish this section with a result on the deterministic Stokes equation, that will be used in the sequel.
Given the deterministic linear  problem
\[\begin{cases}
 \dfrac{dx}{dt}(t)=Ax(t)+f(t), &\qquad t \in (0,T]\\
  x(0)=x_0
 \end{cases}
\]
we represent its mild solution as
\[
x(t)=e^{tA}x_0+\int_0^te^{(t-s)A}f(s)\ ds
\]
and we have ( see \cite{b} Proposition 4.1, based on \cite{dv})
\begin{proposition} \label{isom}
Let $1<p,q,r<\infty$ and $s \in \mathbb R$.\\
For any $f \in L^r(0,T;\mathcal B^{s}_{p\,q})$ and
$x_0 \in \mathcal B^{s+2-\frac{2}{r}}_{p\,r}$,
 there exists a unique solution $x \in W^{1,r}(0,T)\equiv
\{ x \in L^r(0,T;\mathcal B^{s+2}_{p\,q}): \frac{dx}{dt} \in L^r(0,T;\mathcal B^{s}_{p\,q})\}$.

Moreover, the functions $x,\frac{dx}{dt}$ depend continuously on the data
$f$ and $x_0$,
that is there exists a positive constant $C$ such that
$$
\Big( \textstyle \int_0^T ( \|x(t)\|^r_{\mathcal B^{s+2}_{p\,q}}
+  \|\frac{dx}{dt}(t)\|^r_{\mathcal B^{s}_{p\,q}})\,dt\Big)^{1/r}
\leq
\Big( C \textstyle\int_0^T \|f(t)\|^r_{\mathcal B^{s}_{p\,q}} dt\Big)^{1/r} +
\|x_0\|_{\mathcal B^{s+2-\frac{2}{r}}_{p\,r}}
$$
Finally, the space $W^{1,r}(0,T)$
is continuously embedded into the space
$C([0,T];\mathcal B^{s+2-\frac{2}{r}}_{p\,r})$,
that is there exists
a positive constant $C$ such that
\begin{equation}\label{int-cont}
 \|x\|_{C([0,T];\mathcal B^{s+2-\frac{2}{r}}_{p\,r})}  \le \;C\; \|x\|_{W^{1,r}(0,T)}
\end{equation}
and therefore the initial condition makes sense.\\
All the constants  depend only on $p,q,r,s$.
\end{proposition}

%%%%%%%%%%%%%%%%%%%%%%%%%%%%%%%%
\section{The bilinear term}\label{sectionB}
When we study equation  for the auxiliary process $u=v-z$, 
there appears $B(z,z)$. We analyse the space regularity of this term.
Following \cite{AFHK,AC,DPD,AF}, we estimate it with respect to the Gaussian measure $\mu^H$.
\begin{proposition}\label{stimaBz}
Let $\frac 14<H<1$ and
\begin{align}
& \rho<4H-3& \text{ if } \frac 14<H<\frac 12 \label{cond-sotto}\\
                        &      \rho< 2(H-1)& \text{ if } \frac 12\le H<1 \label{cond-sopra}
\end{align}
Then, for any $m \in \mathbb N$ 
\begin{equation}\label{stimaB(z,z)}
\int \|B(z,z)\|^{2m}_{\mathcal H^{\rho}}\ \mu^H(dz)<\infty .
\end{equation}
\end{proposition}
\proof
Let us begin to  perform computations for $m=1$.

First, we explain why we need the lower bound $H>\frac 14$.
By \eqref{componentiBk} we have
\begin{equation}\label{prima-stima}
\begin{split}
\int&\|B(z,z)\|_{\mathcal H^\rho}^2\  \mu^H(dz)
=\sum_{k\in\zo} |k|^{2\rho}\int |B_k(z,z)|^2 \ \mu^H(dz)
\\& 
=\sum_{k\in\zo} |k|^{2\rho} \sum_{h,h^\prime\in\zo; h,h^\prime\neq k}
 \int \gamma_{h,k} z_h z_{k-h}\gamma_{h^\prime,k} \overline{ z_{h^\prime} z_{k-h^\prime}}\ \mu^H(dz)
\\& 
=\sum_{k\in\zo} |k|^{2\rho}\sum_{h \in \zo, h \neq k}
(\gamma^2_{h,k} +\gamma_{h,k}\gamma_{k-h,k}) \int  |z_h|^2 |z_{k-h}|^2 \ \mu^H(dz)
\\& 
=2 \sum_{k\in\zo} |k|^{2\rho}\sum_h \gamma^2_{h,k} \frac {C_H^2}{|h|^{4H} |k-h|^{4H}}
\end{split}
\end{equation}
From \eqref{componentiBk} we have that
$\gamma_{h,k}=\gamma_{k-h,k}$ and $\gamma^2_{h,k}\le |k|^2$; then we can bound 
$\int\|B(z,z)\|_{\mathcal H^\rho}^2\ d\mu^H(z)$ by
\begin{equation}\label{doubles}
\sum_{k\in\zo} |k|^{2\rho+2}\sum_{h\in\zo, h\neq k} \frac 1{|h|^{4H} |k-h|^{4H}}
\end{equation}
For any fixed $k$, the latter series (over $h$) is convergent if  and only if $8H>2$. Therefore we require
\[
H>\frac 14.
\]

The inner series depends on $k$ as proved in Lemma \ref{lemma1} in the
Appendix \ref{ar}. 
Therefore the double series \eqref{doubles} is estimated by
\[
\begin{cases}
\displaystyle\sum_{k\in\zo} |k|^{2\rho+2}\frac 1{ |k|^{8H-2}} &  \text{ if } \frac 14<  H<\frac 12\\
\displaystyle\sum_{k\in\zo} |k|^{2\rho+2}\frac {\ln |k|}{ |k|^{2}} &  \text{ if }  H=\frac 12\\
\displaystyle\sum_{k\in\zo} |k|^{2\rho+2}\frac 1{ |k|^{4H}} &  \text{ if } H> \frac 12
\end{cases}
\]
The first series converges when $\rho<4H-3$, the second one when $\rho<-1$
 and the third one when $\rho<2H-2$. This provides the summability
 \eqref{stimaB(z,z)} under conditions \eqref{cond-sotto}-\eqref{cond-sopra}.

Now, let us consider higher powers $m>1$.
We have that  \eqref{stimaB(z,z)} holds 
also for the other powers, since $\mu^H$ is Gaussian  and therefore 
the higher moments are expressed by means  of the second moments. 
For completeness we provide computations for $m=2$ in Appendix \ref{m2}.
 \hfill\qed%\smartqed
 
Using the stationarity we can write \eqref{stimaB(z,z)} also as
\[
\mathbb E\big[\|B(z(t),z(t))\|^{2m}_{\mathcal H^\rho}\big]=:\tilde e_m<\infty
\]
for an $t\in \mathbb R$.
As an easy consequence,  we obtain
\[
\mathbb E\left[\int_{t_0}^{t_1}\|B(z(t),z(t))\|^{2m}_{\mathcal H^{\rho}} dt \right]
=
(t_1-t_0) \int 
\|B(z,z)\|^{2m}_{\mathcal H^\rho}\ \mu^H(dz)<\infty
\]
for any $\infty<t_0<t_1<\infty$. Hence
%corollary
\begin{corollary}\label{cor-B-inL0T}
Let $m\ge 1$ and  $T>0$. 
Choosing  $\rho$  as in \eqref{cond-sotto}-\eqref{cond-sopra} we get
\[
B(z,z) \in L^m(0,T;\mathcal H^\rho)
\]
$\mathbb P$-a.s.
\end{corollary}
\begin{remark}\label{reg}
Notice that for $\frac 12<H<1$ the quadratic term $B(z,z)$ is in  $L^2(0,T;\mathcal H^{-1})$,
$\mathbb P$-a.s.
\end{remark}

%SECTION
\section{The nonlinear auxiliary equation}\label{sectionNS}
Let $v$ be the unknown for our equation \eqref{NS} 
and let $z$ be the stationary Stokes process given by \eqref{OU}.
The process $u=v-z$ solves the equation
\begin{equation}\label{eq-u}
\frac{du}{dt}=Au-B(u,u)-B(u,z)-B(z,u)-B(z,z).
\end{equation}
For $r<2(H-\frac 12)$ 
we have $z(0)\in \mathcal B^r_{pq}$, $\mathbb P$-a.s. and we take $u(0)=v(0)-z(0)$.

We shall prove that equation \eqref{eq-u} has a local solution when $\frac7{16} <H<\frac 12$
whereas we have a global result when $\frac 12<H<1$.
This implies results for the unknown $v=z+u$.

\subsection{$\frac 14<H<\frac 12$}
We consider  a mild solution $u$ to equation \eqref{eq-u}.
We want to show local existence (and uniqueness) by means of a fixed point argument.
Thus we define the mapping $\mathcal I$
\begin{multline}\label{mappaI}
[\mathcal I(u)](t)=e^{tA}u(0)-
\int_0^t e^{(t-s)A}B(u(s),u(s))\ ds
-\int_0^t e^{(t-s)A}B(z(s),u(s))\ ds\\
-\int_0^t e^{(t-s)A}B(u(s),z(s))\ ds
-\int_0^t e^{(t-s)A}B(z(s),z(s))\ ds
\end{multline}
A  fixed point of $\mathcal I$ is a mild solution of equation \eqref{eq-u}.

Given $T>0$, let 
\[
\mathcal E_T=L^\beta(0,T;\mathcal B^\alpha_{pq})\cap C([0,T];\mathcal B^\sigma_{pq}).
\]
First, we want to show that $\mathcal I: \mathcal E_T\to \mathcal E_T$ for suitable values of the parameters
$\alpha, \beta,\sigma,p,q,H$.

Define
\[
I_0(t)=e^{tA}u_0 .
\]
Given $u_0\in \mathcal B^\sigma_{pq}$, it is an easy result that $I_0 \in \mathcal E_T$ when
\[
\alpha<\sigma+\frac 2\beta .
\]
Indeed, $\|e^{tA}u_0\|_{\mathcal B^\sigma_{pq}}\le \|u_0\|_{\mathcal B^\sigma_{pq}}$; by \eqref{stima-expA} we have
\[
\int_0^T \|e^{tA}u_0\|^\beta_{\mathcal B^\alpha_{pq}}dt\le
C\|u_0\|^\beta_{\mathcal B^\sigma_{pq}} \int_0^T \frac{dt}{t^{\frac{\alpha-\sigma}2\beta}}
\]
and the latter intergal is finite when $\alpha<\sigma+\frac 2\beta$.

To study the integrals involving $B(u,u)$, $B(z,u)$ and $B(u,z)$ we define
\[
I_1(u,\tilde u)(t)=\int_0^t e^{(t-s)A}B(u(s),\tilde u(s))\ ds .
\]
%LEMMA
\begin{lemma}\label{lemma-L}
Let $\alpha,\sigma \in \mathbb R$ and $\beta,p,q\ge 1$
be such that
\[\begin{cases}
\frac 2p+\frac 2\beta<\sigma +1
\\
\alpha<\frac 2p,\quad \sigma<\frac 2p
\\
\alpha+\sigma>0
\end{cases}
\]
If $u \in L^\beta(0,T;\mathcal B^\alpha_{pq})$ and  $\tilde u \in C([0,T];\mathcal B^\sigma_{pq})$, 
then
\[
\|I_1(u,\tilde u)\|_{L^\beta(0,T;\mathcal B^\alpha_{pq})}
\le
CT^{\frac 12-\frac 1p+\frac \sigma 2} \|u\|_{L^\beta(0,T;\mathcal B^\alpha_{pq})} \|\tilde u\|_{C([0,T];\mathcal B^\sigma_{pq})}
\]
and 
\[
\|I_1(\tilde u, u)\|_{L^\beta(0,T;\mathcal B^\alpha_{pq})}
\le
CT^{\frac 12-\frac 1p+\frac \sigma 2} \|u\|_{L^\beta(0,T;\mathcal B^\alpha_{pq})} \|\tilde u\|_{C([0,T];\mathcal B^\sigma_{pq})}
\]
where the constant $C$ is independent of the time $T$.
\end{lemma}
\proof
 We consider the first estimate, since the second one is obtained in the same way interchanging $u$ and $\tilde u$.

We use \eqref{chemin} with $s_1=\alpha$ and $ s_2=\sigma$:
\[
 \|B(u,\tilde u)\|_{\mathcal B^{\alpha+\sigma-\frac2p-1}_{pq}}
\le
C \|u\|_{\mathcal B^\alpha_{pq}}\|\tilde u\|_{\mathcal B^\sigma_{pq}}
\]
where $\alpha<\frac 2p$, $\sigma<\frac 2p$ and $\alpha+\sigma>0$.
Then we get that $B(u,\tilde u)\in L^\beta(0,T;\mathcal B^{\alpha+\sigma-\frac2p-1}_{pq})$. 

Moreover
\[
\|I_1(u,\tilde u)\|^\beta_{L^\beta(0,T;\mathcal B^\alpha_{pq})}
\le
\int_0^T (\int_0^t \|e^{(t-s)A}B(u(s),\tilde u(s))\|_{\mathcal B^\alpha_{pq}}ds)^\beta dt .
\]

Now, we perform estimates using \eqref{stima-expA}   and the H\"older inequality:
\[\begin{split}
\int_0^t &\|e^{(t-s)A}B(u(s),\tilde u(s))\|_{\mathcal B^\alpha_{pq}} ds
\\&\le
\int_0^t \frac C{(t-s)^{\frac12+\frac 1p-\frac \sigma 2}}  \|B(u(s),\tilde u(s))\|_{\mathcal B^{\alpha+\sigma-1-\frac 2p}_{pq}} ds
\\&\le
C\int_0^t \frac 1{(t-s)^{\frac12+\frac 1p-\frac \sigma 2}} 
\|u(s)\|_{\mathcal B^\alpha_{pq}} \|\tilde u(s)\|_{\mathcal B^\sigma_{pq}}ds
\\&\le
C \|\tilde u\|_{C([0,T];\mathcal B^\sigma_{pq})} 
\Big(\int_0^t \frac {ds}{(t-s)^{(\frac12+\frac 1p-\frac \sigma 2)\frac\beta{\beta-1}}}  \Big)^{1-\frac 1\beta}
\Big(\int_0^t \|u(s)\|^\beta_{\mathcal B^\alpha_{pq}}ds\Big)^{\frac 1\beta}
\\&\leq C
t^{\frac 12-\frac 1p+\frac \sigma 2-\frac 1\beta}\|\tilde u\|_{C([0,T];\mathcal B^\sigma_{pq})}  \|u\|_{L^\beta(0,T;\mathcal B^\alpha_{pq})}
\end{split}
\]
Integrating in time over the interval $[0,T]$, we conclude the proof. 
\hfill\qed%\smartqed

Now we consider the other norm for $I_1$.
%LEMMA
\begin{lemma}\label{lemma-C}
Let  $\alpha,\sigma \in \mathbb R$ and $\beta,p,q\ge 1$
be such that
\[
\begin{cases}
\frac 2p+\frac 2q<\alpha+1\\
\beta\ge q\\
\alpha<\frac 2p,\quad \sigma<\frac 2p\\
\alpha+\sigma>0 
\end{cases}
\]
If $u \in L^\beta(0,T;\mathcal B^\alpha_{pq})$ and  $\tilde u \in C([0,T];\mathcal B^\sigma_{pq})$, 
then
\[
\|I_1(u,\tilde u)\|_{C([0,T];\mathcal B^\sigma_{pq})}
\le
C T^{\frac \alpha2+\frac 12-\frac 1p-\frac 1\beta} \|u\|_{L^\beta(0,T;\mathcal B^\alpha_{pq})} \|\tilde u\|_{C([0,T];\mathcal B^\sigma_{pq})}
\]
and
\[
\|I_1(\tilde u,u)\|_{C([0,T];\mathcal B^\sigma_{pq})}
\le
C T^{\frac \alpha2+\frac 12-\frac 1p-\frac 1\beta} \|u\|_{L^\beta(0,T;\mathcal B^\alpha_{pq})} \|\tilde u\|_{C([0,T];\mathcal B^\sigma_{pq})}
\]
where  the constant $C$ is independent of the time $T$.
\end{lemma}
\proof
First, from the previous  proof we know 
that $B(u,\tilde u)\in L^\beta(0,T; \mathcal B^{\alpha+\sigma-\frac 2p-1}_{pq})$; when $\beta\ge q$ we also have 
$B(u,\tilde u)\in L^q(0,T; \mathcal B^{\alpha+\sigma-\frac 2p-1}_{pq})$ and Proposition \ref{isom}
provides $I_1(u,\tilde u)\in C([0,T];\mathcal B^{\alpha+\sigma-\frac 2p+1-\frac 2q}_{pq}) $ and finally we use that
$\mathcal B^{\alpha+\sigma-\frac 2p+1-\frac 2q}_{pq}\subseteq \mathcal B^{\sigma}_{pq}$ when 
$\frac 2p+\frac 2q\le \alpha+1 $.

Now, we perform estimates using
\eqref{stima-expA} and the H\"older inequality.
\[\begin{split}
\|I_1(u,\tilde u)(t)\|_{\mathcal B^\sigma_{pq}}&\le
\int_0^t \|e^{(t-s)A}B(u(s),\tilde u(s))\|_{\mathcal B^\sigma_{pq}}ds
\\&\le
\int_0^t \frac C{(t-s)^{\frac {1}2+\frac 1p-\frac \alpha 2}} \|B(u(s),\tilde u(s))\|_{\mathcal B^{\alpha+\sigma-\frac2p-1}_{pq}}ds
\\&
\le\int_0^t \frac C{(t-s)^{\frac {1}2+\frac 1p-\frac \alpha 2}}\|u(s)\|_{\mathcal B^\alpha_{pq}}\|\tilde u(s)\|_{\mathcal B^\sigma_{pq}} ds
\\&
\le \|\tilde u\|_{C([0,T];\mathcal B^\sigma_{pq})} \int_0^t \frac C{(t-s)^{\frac {1}2+\frac 1p-\frac \alpha 2}}\|u(s)\|_{\mathcal B^\alpha_{pq}}ds
\\&\le
C \|\tilde u\|_{C([0,T];\mathcal B^\sigma_{pq})} \|u\|_{L^\beta(0,T;\mathcal B^\alpha_{pq})} 
(\int_0^t \frac{ds}{(t-s)^{\frac \beta{\beta-1}(\frac {1}2+\frac 1p-\frac \alpha 2)}})^{1-\frac 1\beta}
\end{split}
\]
The latter intergral is finite when $(\frac{\beta}{\beta-1})(\frac {1}2+\frac 1p-\frac \alpha 2)<1$, 
i.e.
\begin{equation}\label{cond-Linfty}
\frac 2\beta+\frac 2p<\alpha+1
\end{equation}
This inequality is true when $\beta\ge q $ and $\frac 2q+\frac 2p< \alpha+1$, 
that is our assumptions imply \eqref{cond-Linfty}.

Computing the time interval and taking the supremum over 
$t\in [0,T]$ we get the required estimate. 
\hfill\qed%\smartqed

For the integral involving $B(z,z)$ we define the process 
\[
I_2(t)=\int_0^t e^{(t-s)A}B(z(s),z(s))ds,\qquad t\ge 0
\]
where $z$ is the Stokes process given in \eqref{OU}.
%LEMMA
\begin{lemma}\label{lemma-B}
Let $\frac 14<H<\frac 12$, $\beta, p\ge 1$, $q\ge 2$
and $\alpha, \sigma \in \mathbb R$ be such that 
\[
\alpha\le \sigma+1 \ , \qquad \sigma<4H-2.
\]
Then 
$I_2\in \mathcal E_T$, $\mathbb P$-a.s.
\end{lemma}
\proof We proceed pathwise. First we show that $I_2\in L^\beta(0,T;\mathcal B^\alpha_{pq})$.
From Corollary \ref{cor-B-inL0T} we know that the paths of $B(z,z)$ are in $L^\beta(0,T;\mathcal B^{\sigma-1}_{pq})$ for any $\beta\ge 1$ and for 
$\sigma<4H-2$.
Therefore, according to Proposition \ref{isom}
the paths of $I_2$ are in $L^\beta(0,T;\mathcal B^{\sigma+1}_{pq})$. When $\alpha\le \sigma+1 $, the embedding theorem 
gives $I_2 \in L^\beta(0,T;\mathcal B^{\alpha}_{pq})$.

Now, we show that 
$I_2 \in C([0,T];\mathcal B^\sigma_{pq})$.
Again by Corollary \ref{cor-B-inL0T}, for any $q\ge 1$ and $\sigma<4H-2$
the paths of $B(z,z)$ are  in $L^q(0,T;\mathcal B^{\sigma-1}_{pq})$. 
We bear in mind Proposition \ref{isom} and we get
that
$I_2\in C([0,T];\mathcal B^{\sigma+1-\frac2q}_{pq})$. When $q\ge 2$ this finishes the proof.  
\hfill\qed%\smartqed

Summing up, we have proved estimates for all the terms in the r.h.s. of \eqref{mappaI}.
Let us point out that merging these results 
we have to satisfy two conditions:
\[
\sigma<2(H-\frac12)
\] 
which comes from Proposition \ref{pro-zLm} and provides
$z \in C([0,T];\mathcal H^\sigma)$, $\mathbb P$-a.s., so to apply Lemma \ref{lemma-L} and \ref{lemma-C} for 
the integrals involving $B(z,u)$ and $B(u,z)$, and
\[
\sigma <4H-2
\]
which comes from Lemma \ref{lemma-B} to estimate the integral involving $B(z,z)$.

When $H<\frac 12$, the latter condition is stronger and we will write only this one in the following.

\begin{proposition}
Let $\frac 14<H<\frac12$, $\alpha,\sigma \in \mathbb R$, $\beta,p\ge 1$ and $q\ge 2$ 
be such that
\begin{align}
%\begin{cases}
\frac 2p+\frac 2\beta&<\sigma +1  \label{uno}
\\
\frac 2p+\frac 2q&<\alpha+1 \label{due}
\\
\beta&\ge q \label{tre}
\\
\alpha&<\frac 2p \label{quattro}
\\
\sigma&<\frac 2p    \label{cinque}
\\
\alpha+\sigma&>0 \label{sei}
\\
\alpha&<\sigma+\frac 2\beta \label{sette}
\\
\alpha&\le \sigma+1 \label{otto}
\\
 \sigma&<4(H-\frac 12) \label{nove}
%\end{cases}
\end{align}
Then, for any finite $T$ we have that $\mathcal I: \mathcal E_T\to \mathcal E_T$, $\mathbb P$-a.s..
\end{proposition}

\begin{remark}[How to fulfil conditions]
Notice that by condition \eqref{nove} we have $\sigma<0$ when $H<\frac 12$.  Therefore, 
condition \eqref{sei} requires $\alpha>0$. 

Moreover, conditions \eqref{sei} and \eqref{otto} provide
\[
-\sigma<\alpha\le \sigma+1;
\]
thus it is necessary that $\sigma>-\frac12$, i.e. $H>\frac38$.

Actually we are going to show that $H$ must be bigger than $\frac 38$ in order to satisfy all the conditions \eqref{uno}-\eqref{nove}.
Indeed, we write a system equivalent to the previous one. Condition \eqref{cinque} is trivially satisfied when $\sigma<0$ and can be neglected. Taking $\beta=q$, condition \eqref{due} is weaker than condition \eqref{uno}, that is \eqref{uno} implies \eqref{due}.
In addition, since condition \eqref{uno} requires $\frac 2\beta<1$, we have that condition \eqref{otto} is weaker than condition \eqref{sette}, that is \eqref{sette} implies \eqref{otto}. Therefore, in the case  $\beta=q$ the previous system of conditions is equivalent to
\begin{align}
\frac 2p+\frac 2\beta&<\sigma +1  \tag{\ref{uno}}
\\
\beta&= q \tag{\ref{tre}'}
\\
\alpha&<\frac 2p \tag{\ref{quattro}}
\\
\alpha+\sigma&>0  \tag{\ref{sei}}
\\
\alpha&<\sigma+\frac 2\beta  \tag{\ref{sette}}
\\
 \sigma&<4(H-\frac 12)   \tag{\ref{nove}}
\end{align}
which is simpler to analyse.  Let us notice that
\[
-3\sigma\underset{\text{by }\eqref{sei}}{<} 
\alpha+ \alpha-\sigma\underset{\text{by }\eqref{quattro} \text{ and } \eqref{sette}}{<}
 \frac 2p+\frac2 \beta\underset{\text{by }\eqref{uno}}{<}
  \sigma+1
\]
This sequence of inequalities is meaningful only when $-3\sigma< \sigma+1$, i.e. $\sigma>-\frac 14$.
Taking into account the last condition \eqref{nove}, we see that in order to fulfil all the above conditions it is necessary that $H>\frac7{16}$.

Therefore it is possible to fulfil all the conditions when $\frac 7{16}<H<\frac 12$. Setting $H=\frac 7{16}+c$ with 
$0<c<\frac 1{16}$, we can choose for instance
\[
\alpha=\frac 14-2c,\qquad 
\sigma=-\frac 14+3c,\qquad
\frac 2\beta=\frac 2q=\frac 12-4c,\qquad
\frac 2p=\frac 14-c
\]
in order to satisfy \eqref{uno}-\eqref{nove}.
\end{remark}

Now we can prove the local existence result for $u$, proving
that  $\mathcal I$ is a contraction for $T$ small enough.
\begin{proposition}
Let  $\frac 7{16}< H<\frac 12$ and 
 the parameters fulfil  the conditions \eqref{uno}-\eqref{nove}.
Then, given $u_0 \in \mathcal B^\sigma_{pq}$
there exist a  stopping time $\tau\in ]0,T]$ and for $\mathbb P$-a.e. $\omega$ a unique mild solution $u(\omega,\cdot)$ 
 of equation \eqref{eq-u} with values in 
$C([0,\tau(\omega)];\mathcal B^\sigma_{pq})\cap L^\beta(0,\tau(\omega);\mathcal B^\alpha_{pq})$.
\end{proposition}
\proof
Using the bilinearity of the operator $B$, we get 
\[\begin{split}
\mathcal I(u_1)(t)-\mathcal I(u_2)(t)=&-\int_0^t e^{(t-s)A}B(u_1(s),u_1(s)-u_2(s))ds\\
&-\int_0^t e^{(t-s)A}B(u_1(s)-u_2(s),u_2(s))ds
\\&-\int_0^t e^{(t-s)A}B(z(s),u_1(s)-u_2(s))ds\\
&-\int_0^t e^{(t-s)A}B(u_1(s)-u_2(s),z(s))ds
\end{split}\]
Let us work in the subspace of $\mathcal E_T$ with $\|u\|_{\mathcal E_T}\le M$.
The initial data $u(0) \in \mathcal B^\sigma_{pq}$ is fixed.
Therefore, according to  Lemma \ref{lemma-L} and  Lemma \ref{lemma-C}   we have
\begin{multline}
\|\mathcal I(u_1)-\mathcal I(u_2)\|_{\mathcal E_T}
\\\le 
\overline C  \big(T^{\frac 12-\frac 1p+\frac \sigma 2}+T^{\frac \alpha2+\frac 12-\frac 1p-\frac 1\beta}\big)
(M+\|z\|_{C([0,T];\mathcal B^\sigma_{pq})})
\|u_1-u_2\|_{\mathcal E_T}
\end{multline}
for a suitable constant $\overline C$ independent of $T$.

When $T$ is such that
\begin{equation}\label{ineqT}
\overline C  \big(T^{\frac 12-\frac 1p+\frac \sigma 2}+T^{\frac \alpha2+\frac 12-\frac 1p-\frac 1\beta}\big)
(M+\|z\|_{C([0,T];\mathcal B^\sigma_{pq})})<1
\end{equation}
the mapping $\mathcal I$ is a contraction and hence has a unique fixed point, which is the unique solution 
of equation \eqref{eq-u}. 

Notice that $T$ is a random time, since inequality \eqref{ineqT} involves the random process $z$. It can be chosen to be a stopping time. 
\hfill\qed%\smartqed

Since $v=z+u$, we also get existence of a local mild solution $v$ to equation 
\eqref{NS} where the bilinear term $B(v,v)$ has to be understood as the sum of four terms, that is
\begin{multline}\label{nuovoNS}
dv(t)-Av(t)\ dt= -B(u(t),u(t))\ dt-B(u(t),z(t))\ dt
\\-B(z(t),u(t))\ dt-B(z(t),z(t))\ dt +dw^H(t)
\end{multline}
\begin{theorem} 
Let  $\frac 7{16}< H<\frac 12$ and 
 the parameters fulfil  the conditions \eqref{uno}-\eqref{nove}.
Then, given $v_0 \in \mathcal B^\sigma_{pq}$
there exist a  stopping time $\tau\in ]0,T]$ and 
 for $\mathbb P$-a.e. $\omega$ a unique mild solution $v(\omega,\cdot)$ 
 of equation \eqref{nuovoNS} with values in 
$C([0,\tau(\omega)];\mathcal B^\sigma_{pq})$.
\end{theorem}
\begin{remark} 
We cannot get a global existence result for $v$ as in \cite{AC,DPD}; indeed, when $H=\frac 12$ the
Gaussian  measure
$\mu^H$ defined by \eqref{measure}
 is invariant for the Navier-Stokes equation \eqref{NS}. 
This allows  to define $B(v,v)$ and to get  global existence.
However, when $H\neq \frac12$  the
 measure
$\mu^H$ is invariant for the Stokes equation \eqref{eq-lin} but not 
 for the Navier-Stokes equation \eqref{NS}; this depends eventually on the fact that for $0<H<1$ 
the Gaussian measure $\mu^H$ is (formally) invariant for the deterministic Euler dynamics 
\[
\frac{dv}{dt}= -B(v,v)
\]
only when $H=\frac 12$ (and in this case $\mu^{\frac 12}$ is called the enstrophy measure, see \cite{AFHK}). 
\end{remark}

%%%%%%%%%%%%%%%%%%%%%%%%%%%%%%%%%%
\subsection{$\frac12<H<1$}
When $\frac12<H<1$ the fBm $w^H$ and the Stokes process $z$ are more regular and 
we expect more regularity of the processes $u$ and $v$ too.
Actually we can obtain an a priori energy estimate; this will lead to global existence.
Let us notice that now we deal with  solutions which are weak in the sense of PDE's; for instance the solution $u$ of equation \eqref{eq-u} has paths at least in  $L^\infty(0,T;\mathcal H^0)\cap L^2(0,T;\mathcal H^{1})$
and  fulfils for any $t>0$ and any $\varphi\in \mathcal H^1$
\begin{multline*}
\langle u(t)-u(0),\varphi\rangle-\int_0^t \langle Au(s), \varphi\rangle ds +\int_0^t
\langle B(u(s),\varphi), u(s)\rangle ds
\\+\int_0^t \langle B(u(s),\varphi) ,z(s)\rangle ds
+\int_0^t \langle B(z(s),\varphi), u(s)\rangle ds
=-\int_0^t \langle B(z(s),z(s)),\varphi\rangle ds
\end{multline*}
$\mathbb P$-a.s.. This is obtained from \eqref{eq-u} by using \eqref{tril}.

Since  the paths of the process $u$ are in $L^\infty(0,T;\mathcal H^0)\cap L^2(0,T;\mathcal H^{1})$ and those of $z$ are in 
$C([0,T];\mathcal H^\sigma)$ for some $\sigma>0$,
then all the terms in the latter relationship are well defined. Let us check the trilinear terms, by using H\"older inequality, interpolation inequality and Sobolev embeddings:
\[
|\langle B(u(s),\varphi), u(s)\rangle|
\le \|u(s)\|_{L^4}^2\|\varphi\|_{\mathcal H^1}
\le C\|u(s)\|_{\mathcal H^0}  \|u(s)\|_{\mathcal H^1} \|\varphi\|_{\mathcal H^1}
\]
\[
\begin{split}
|\langle B(u(s),\varphi), z(s)\rangle|
&\le \|u(s)\|_{L^{\frac 2\sigma}}\|\varphi\|_{\mathcal H^1} \|z(s)\|_{L^{\frac 2{1-\sigma}}}\; \text{ for } 0<\sigma<1
\\&\le C\|u(s)\|_{\mathcal H^1}   \|\varphi\|_{\mathcal H^1}\|z(s)\|_{\mathcal H^\sigma}
\end{split}
\]
The third trilinear term can be dealt with as with the second term. And finally the latter term is well defined as soon as 
$B(z(s),z(s))\in L^1(0,T;\mathcal H^{-1})$ (see Remark \ref{reg}).

Now, let $0<\sigma<2(H-\frac 12)$ for $\frac 12<H<1$.
Taking the $\mathcal H^0$-scalar product of  equation \eqref{eq-u} with  $u$, we get the usual 
energy estimate (see \cite{Temam}). We make use of \eqref{Bnull} and \eqref{stimeGiga}:
\[\begin{split}
\frac 12 \frac{d}{dt} &\|u(t)\|_{\mathcal H^0}^2+\|\nabla u(t)\|^2_{L^2}
\\&
=-\langle B(u(t)+z(t),u(t)),u(t)\rangle
-\langle B(u(t),z(t)),u(t)\rangle-\langle B(z(t),z(t)),u(t)\rangle
\\
&=-\langle B(u(t), z(t)),u(t)\rangle-\langle B(z(t),z(t)),u(t)\rangle
\\
&\le
\|B(u(t), z(t))\|_{\mathcal H^{-1}} \|u(t)\|_{\mathcal H^1}
+\|B(z(t), z(t))\|_{\mathcal H^{-1}} \|u(t)\|_{\mathcal H^1}
\\
&\le C\|u(t)\|_{\mathcal H^{1-\sigma}}  \|z(t)\|_{\mathcal H^\sigma}\|u(t)\|_{\mathcal H^1}
+\frac 14 \|u(t)\|_{\mathcal H^1}^2 +C  \|B(z(t),z(t))\|_{\mathcal H^{-1}}^2
\end{split}
\]
Moreover, by interpolation and Young inequality
\[\begin{split}
\|u\|_{\mathcal H^{1-\sigma}}  \|z\|_{\mathcal H^\sigma}\|u\|_{\mathcal H^1}
&\le
C\|u\|^\sigma_{\mathcal H^0} \|u\|^{1-\sigma}_{\mathcal H^1}  \|z\|_{\mathcal H^\sigma}\|u\|_{\mathcal H^1}
\\&=
C\|u\|^\sigma_{\mathcal H^0} \|u\|^{2-\sigma}_{\mathcal H^1}  \|z\|_{\mathcal H^\sigma}
\\&\le
\frac 14 \|u\|^{2}_{\mathcal H^1} +C\|u\|^2_{\mathcal H^0}\|z\|^{\frac 2\sigma}_{\mathcal H^\sigma}
\end{split}\]
Since $\|u\|^2_{\mathcal H^1}=\|u\|^2_{\mathcal H^0}+\|\nabla u\|^2_{L^2}$, collecting all the estimates 
we have found
\[
\frac{d}{dt} \|u(t)\|_{\mathcal H^0}^2+\|\nabla u(t)\|^2_{L^2}
\le C(1+\|z\|^{\frac 2\sigma}_{\mathcal H^\sigma})\|u\|^2_{\mathcal H^0} +
C  \|B(z(t),z(t))\|_{\mathcal H^{-1}}^2
\]
According to Remark \ref{reg}, $B(z,z)\in L^2(0,T;\mathcal H^{-1})$.
Moreover $z \in C([0,T];\mathcal H^\sigma)$ by Proposition \ref{pro-zLm}.
This provides as usual by means of Gronwall Lemma
that $u \in L^\infty(0,T;\mathcal H^0)\cap L^2(0,T;\mathcal H^1)$, $\mathbb P$-a.s.
The reader can see all the details of this standard procedure 
in \cite{Temam}. First one has to work on the finite dimensional approximation and then pass to the limit.
By interpolation
\begin{equation}\label{u-reg}
u \in L^{\frac 2\sigma}(0,T;\mathcal H^\sigma).
\end{equation}
Hence, $v=z+u \in L^\infty(0,T;\mathcal H^0)\cap L^{\frac 2\sigma}(0,T;\mathcal H^\sigma)$.

We can improve the estimates, now getting
$u \in C([0,T];\mathcal H^\sigma)\cap L^2(0,T;\mathcal H^{1+\sigma})$. 
This gives global existence for the process $v=z+u$ in the space $C([0,T];\mathcal H^\sigma)$ for
$0<\sigma<2(H-\frac 12)$. Therefore the term $B(v,v)$ is well defined. Actually 
the process $v$ is a weak solution (in the sense of PDE's) to equation \eqref{NS}, that is it solves
for any $t>0$ and any $\varphi\in \mathcal H^{2-\sigma}$
\begin{multline*}
\langle v(t)-v_0,\varphi\rangle-\int_0^t \langle Av(s), \varphi\rangle ds +\int_0^t
\langle B(v(s),\varphi), v(s)\rangle ds
=\langle w^H(t),\varphi\rangle 
\end{multline*}
$\mathbb P$-a.s.. We leave to the reader to check that all term are well defined (use that $2-\sigma>1$).
 Similarly, the process $z$  can be considered as a weak solution of the 
stochastic Stokes equation \eqref{eq-lin}.

\begin{theorem}[Global existence] \label{th-global-ex}
Let $\frac 12<H<1$ and 
\begin{equation}
0<\sigma<2(H-\frac 12).
\end{equation}
Given $v_0\in\mathcal H^\sigma$ there exists a 
$C([0,T];\mathcal H^\sigma)\cap L^2(0,T;\mathcal H^{1+\sigma})$-valued process $u$ 
solving equation \eqref{eq-u} with $u(0)=v_0-z(0)$.
Therefore there exists a $C([0,T];\mathcal H^\sigma)$-valued process $v$ solving 
equation \eqref{NS} with $v(0)=v_0$.
\end{theorem}
\proof
We have to work on
\[\begin{split}
\frac 12 \frac{d}{dt} &\|u(t)\|_{\mathcal H^\sigma}^2+\|\nabla u(t)\|^2_{H^\sigma}
\\&
=-\langle B(u(t)+z(t),u(t)),(-A)^{\sigma}u(t)\rangle
-\langle B(u(t),z(t)),(-A)^{\sigma}u(t)\rangle
\\&\qquad-\langle B(z(t),z(t)),(-A)^{\sigma}u(t)\rangle
\end{split}
\]
We performe the estimates on the terms in the r.h.s.. Notice that $0<\sigma<1$.

From \eqref{stimeGiga}, the interpolation inequality 
$\|u\|_{\mathcal H^{1}}\le C\|u\|^{\sigma}_{\mathcal H^{\sigma}}\|u\|^{1-\sigma}_{\mathcal H^{1+\sigma}}$
and Young inequality we get
\[\begin{split}
\langle B(u,u),(-A)^{\sigma}u\rangle
&\le
\|B(u,u)\|_{\mathcal H^{\sigma-1}} \|u\|_{\mathcal H^{\sigma+1}}
\\&
\le C \|u\|_{\mathcal H^{\sigma}}\|u\|_{\mathcal H^{1}}\|u\|_{\mathcal H^{1+\sigma}}
\\&
\le C \|u\|^{1+\sigma}_{\mathcal H^{\sigma}}\|u\|^{2-\sigma}_{\mathcal H^{1+\sigma}}
\\&\le
\frac 18 \|u\|^{2}_{\mathcal H^{1+\sigma}}+C \|u\|^{\frac 2 \sigma}_{\mathcal H^{\sigma}}\|u\|^2_{\mathcal H^{\sigma}}
\end{split}\]

\[\begin{split}
\langle B(z,u),(-A)^{\sigma}u\rangle
&\le
\|B(z,u)\|_{\mathcal H^{\sigma-1}} \|u\|_{\mathcal H^{\sigma+1}}
\\&
\le C \|z\|_{\mathcal H^{\sigma}}\|u\|_{\mathcal H^{1}}\|u\|_{\mathcal H^{1+\sigma}}
\\&\le
C  \|z\|_{\mathcal H^{\sigma}}\|u\|^{\sigma}_{\mathcal H^{\sigma}}\|u\|^{2-\sigma}_{\mathcal H^{1+\sigma}}
\\&\le 
\frac 18 \|u\|^2_{\mathcal H^{1+\sigma}}+ C \|z\|_{\mathcal H^{\sigma}}^{\frac 2 \sigma}\|u\|^2_{\mathcal H^{\sigma}}
\end{split}\]

\[\begin{split}
\langle B(u,z),(-A)^{\sigma}u\rangle
&\le
\|B(u,z)\|_{\mathcal H^{\sigma-1}} \|u\|_{\mathcal H^{\sigma+1}}
\\&
\le C \|u\|_{\mathcal H^{1}}\|z\|_{\mathcal H^{\sigma+\epsilon}}\|u\|_{\mathcal H^{1+\sigma}}
\\&\le C \|u\|^{\sigma}_{\mathcal H^{\sigma}}
\|z\|_{\mathcal H^{\sigma+\epsilon}}\|u\|^{2-\sigma}_{\mathcal H^{1+\sigma}}
\\&\le
\frac 18  \|u\|^2_{\mathcal H^{1+\sigma}}+ C \|z\|_{\mathcal H^{\sigma+\epsilon}}^{\frac 2 \sigma}\|u\|^2_{\mathcal H^{\sigma}}
\end{split}\]
where $0<\epsilon\ll 1$, and 
\[\begin{split}
\langle B(z,z),(-A)^{\sigma}u\rangle
&\le
\|B(z,z)\|_{\mathcal H^{\sigma-1}} \|u\|_{\mathcal H^{\sigma+1}}
\\&
\le \frac 18  \|u\|^2_{\mathcal H^{1+\sigma}}+ C \|B(z,z)\|^2_{\mathcal H^{\sigma-1}}
\end{split}\]

Now, $ \|u\|^2_{\mathcal H^{1+\sigma}}=\|u\|^2_{\mathcal H^\sigma}+\|\nabla u\|^2_{H^\sigma}$.
Summing up
\[
 \frac{d}{dt} \|u(t)\|_{\mathcal H^\sigma}^2+\|\nabla u(t)\|^2_{H^\sigma}
\le C
(1+\|u\|^{\frac 2 \sigma}_{\mathcal H^{\sigma}}+\|z\|_{\mathcal H^{\sigma+\epsilon}}^{\frac 2 \sigma}
+\|B(z,z)\|^2_{\mathcal H^{\sigma-1}}) \|u\|^2_{\mathcal H^{\sigma}} .
\]

Now we bear in mind \eqref{u-reg},  Proposition \ref{pro-zLm} with $\epsilon\ll 1$ and
Proposition \ref{stimaBz} to deal with the sum in the r.h.s..
So by means of  Gronwall lemma 
we conclude that
$u \in L^\infty(0,T;\mathcal H^{\sigma})\cap L^2(0,T;\mathcal H^{\sigma+1})$. 
In addition, by means of the previous estimates we get
\[
\frac{du}{dt}=Au-B(u,u)-B(u,z)-B(z,u)-B(z,z)
\in L^2(0,T;\mathcal H^{\sigma-1}).
\]
Since $u\in  L^2(0,T;\mathcal H^{\sigma+1})$, one gets that 
$u \in C([0,T];\mathcal H^{\sigma})$ (see \cite{Temam}). Hence
$v=u+z\in  C([0,T];\mathcal H^{\sigma})$. 
\hfill\qed%\smartqed

\medskip
The solution obtained is also unique. We have a pathwise uniqueness result.
\begin{theorem}[Uniqueness]
Let $\frac 12<H<1$ and $0<\sigma<2(H-\frac 12)$.
Given $v_0\in\mathcal H^\sigma$ there exists a unique 
 $C([0,T];\mathcal H^\sigma)$-valued process solving \eqref{NS}.
\end{theorem}
\proof
Let $v_1,v_2 \in C([0,T];\mathcal H^\sigma)$ be solutions of \eqref{NS}.
Then the difference $V=v_1-v_2$ fulfils
\begin{equation}\label{eq-V}
\frac {dV}{dt}+AV=-B(v_1,v_1)+B(v_2,v_2)
\end{equation}
with $V(0)=0$. 
We are going to prove that $V(t)=0$ for all $t\ge0$ and this is obtained by means of the a priori estimate
of the energy. Actually the paths of $V$ are more regular than those of $v_1$ and $v_2$, since the
noise term has desappeared in \eqref{eq-V}; this was remarked
already in \cite{Fe03}.
More precisely, we state that any solution $V$ of
\eqref{eq-V} with $V(0)=0$  is such that 
$ V\in C([0,T];\mathcal H^0)\cap L^2(0,T; \mathcal H^1)$ 
and  $\frac {dV}{dt}\in L^2(0,T;\mathcal H^{-1})$;
therefore the equality $\frac{d}{dt}\|V(t)\|^2_{\mathcal H^0}=2\langle \frac
{dV}{dt}(t),V(t)\rangle$ holds  and the energy estimates (coming later) are justified.

Indeed we are given $v_1,v_2 \in C([0,T];\mathcal H^\sigma)$ with 
$\sigma \in (0,1)$.
The r.h.s.  of \eqref{eq-V} 
belongs to $L^\infty(0,T;\mathcal H^{2\sigma-2})$, since 
\[
\|B(v_i,v_i)\|_{\mathcal H^{2\sigma-2}} \le C \|v_i\|_{\mathcal H^\sigma}^2
\]
by \eqref{stimeGiga}.
According to  Proposition \ref{isom} (used with vanishing initial data and $r=2$) we get that any solution $V$ will be in 
$L^2(0,T;\mathcal H^{2\sigma})$. 

If $2\sigma\ge 1$
(i.e. when $\frac 12\le \sigma<1$) we have obtained that $V \in L^2(0,T; \mathcal
H^1)$; moreover,  $\frac {dV}{dt}=-AV+B(v_1,v_1)-B(v_2,v_2)$ and
therefore
 $\frac {dV}{dt}\in L^2(0,T;\mathcal H^{-1})$.

Otherwise, when $0<\sigma<\frac 12$ we proceed as follows; 
by the bilinearity of $B$ we get that  \eqref{eq-V} can
be written as
\begin{equation}\label{eq-diffV}
\frac {dV}{dt}+AV=-B(v_1,V)-B(V,v_2)
\end{equation}
Let us look at the regularity of the r.h.s., knowing that 
$v_1,v_2 \in C([0,T];\mathcal H^\sigma)$ and 
$V\in C([0,T];\mathcal H^\sigma)\cap L^2(0,T;\mathcal H^{2\sigma})$.

Thanks to \eqref{stimeGiga} we get, for $0<\sigma<\frac 12$
\[
\|B(v_1,V)\|_{\mathcal H^{3\sigma-2}}\le C \|v_1\|_{\mathcal H^{\sigma}} \|V\|_{\mathcal H^{2\sigma}}
\]
\[
\|B(V,v_2)\|_{\mathcal H^{3\sigma-2}}\le C  \|V\|_{\mathcal H^{2\sigma}} \|v_2\|_{\mathcal H^{\sigma}} 
\]
Hence the r.h.s. of  \eqref{eq-diffV} belongs to $L^2(0,T;\mathcal
H^{3\sigma-2})$ and therefore thanks to Proposition \ref{isom} any
solution $V$  belongs to  $L^2(0,T;\mathcal H^{3\sigma})$. 

If $3\sigma\ge 1$ (i.e. when $\frac 13\le \sigma<\frac 12$) we have
obtained that  $V \in L^2(0,T; \mathcal H^1)$ and moreover the
r.h.s. of \eqref{eq-diffV} belongs to 
$L^2(0,T; \mathcal H^{3\sigma-2})\subseteq  L^2(0,T; \mathcal H^{-1})$. Hence
we conclude as in the previous case about $\frac {dV}{dt}$.

Otherwise, for
smaller values of $\sigma$ we proceed again with the bootstrap
argument.
We conclude  that, given
 $\sigma \in (0,1)$ and $v_i \in C([0,T]; \mathcal H^\sigma)$, any solution $V$ to
\eqref{eq-diffV} is in $ L^2(0,T; \mathcal H^1)$ and  $\frac {dV}{dt}\in L^2(0,T;\mathcal H^{-1})$.

Now, we look for the a priori energy estimate.
Keeping in mind \eqref{Bnull}, \eqref{stimeGiga} and the interpolation inequality
$\|V\|_{\mathcal H^{1-\sigma}}\le C  \|V(t)\|^\sigma_{\mathcal H^{0}} \|V(t)\|^{1-\sigma}_{\mathcal H^1}$, we get
\[\begin{split}
\frac 12 \frac{d}{dt}\|V(t)\|^2_{\mathcal H^0}+\|\nabla V(t)\|^2_{L^2}
&= -
\langle B(v_1(t),V(t)), V(t) \rangle-
\langle B(V(t),v_2(t)),  V(t) \rangle
\\
&=\langle B(V(t),V(t)),v_2(t) \rangle\\
&\le \|B(V(t),V(t))\|_{\mathcal H^{-\sigma}} \|v_2(t)\|_{\mathcal H^{\sigma}}\\
&\le C \|V(t)\|_{\mathcal H^{1-\sigma}} \|V(t)\|_{\mathcal H^1} \|v_2(t)\|_{\mathcal H^{\sigma}}\\
&\le C  \|V(t)\|^\sigma_{\mathcal H^{0}} \|V(t)\|^{2-\sigma}_{\mathcal H^1} \|v_2(t)\|_{\mathcal H^{\sigma}}\\
&\le
\frac 12\|V(t)\|^{2}_{\mathcal H^1}+C \|v_2(t)\|_{\mathcal H^{\sigma}}^{\frac 2\sigma}
\|V(t)\|^2_{\mathcal H^{0}}\\
&=\frac 12\|V(t)\|^{2}_{\mathcal H^0}+\frac 12\|\nabla V(t)\|^{2}_{L^2}+
C \|v_2(t)\|_{\mathcal H^{\sigma}}^{\frac 2\sigma}\|V(t)\|^2_{\mathcal H^{0}}
\end{split}
\]
From
\[
\frac{d}{dt}\|V(t)\|^2_{\mathcal H^0}\le 
C(1+ \|v_2(t)\|_{\mathcal H^{\sigma}}^{\frac 2\sigma})\|V(t)\|^2_{\mathcal H^{0}}
\]
we conclude by Gronwall lemma that $\displaystyle{\sup_{0\le t\le T}}
\|V(t)\|_{\mathcal H^0}=0$. 
This proves pathwise uniqueness. \hfill\qed%\smartqed

\appendix
\section{Auxiliary results}\label{ar}
We prove some results about convergence of series.
%1
\begin{lemma}\label{lemma1}
For any $k \in \zo$,  the series
\[
\sum_{h\in\zo, h\neq k} \frac 1{|h|^{4H} |k-h|^{4H}}
\]
converges if $H>\frac 14 $ 
and its sum $S_1(k)$ depends on $k$ as follows
\[
S_1(k)\le \begin{cases}
M_H\frac 1{ |k|^{8H-2}} &  \text{ if } \frac 14<  H<\frac 12\\
M\frac {\ln |k|}{ |k|^{2}} &  \text{ if }  H=\frac 12\\
M_H\frac 1{ |k|^{4H}} &  \text{ if } H> \frac 12
\end{cases}
\]
for some positive constants $M$ and $M_H$ independent of $k$.
\end{lemma}
\proof
The series  can be estimated by the following integral
\[
\iint_{D_k} \frac {dx\ dy}{[x^{2}+ y^{2}]^{2H} \ [(x-k^{(1)})^{2}+ (y-k^{(2)})^2]^{2H}}
\]
over the region $D_k$ which is the plane without two unitary balls around the points $k$ and $0$. Therefore it is enough to evaluate
\[
I_1(a)=\iint_{\tilde D_a} \frac {dx\ dy}{[x^{2}+ y^{2}]^{2H} \ [(x-a)^{2}+ y^2]^{2H}}
\]
for $a\ge 2$, where $\tilde D_a=\mathbb R^2\setminus(B_1((0,0))\cup B_1((a,0)))$ and 
$B_r((c,d))$  denotes the ball of radius $r$ with center $(c,d)$.

Now we make a  change of variables: 
$u=\frac xa$ and $v=\frac ya$. The domain $D_a$ becomes the domain 
$R_a=\mathbb R^2\setminus (B_{\frac 1a}((0,0))\cup B_{\frac 1a}((1,0)))$.  
Hence
\[
I_1(a)=
\frac{1}{a^{8H}} \iint_{R_a}
 \frac{a^2}{(u^{2}+ v^{2})^{2H} \ ((u-1)^{2}+ v^{2})^{2H}} \ du \ dv .
 \]
We split the integral region into three disjoint  regions: 
$R^0_a=B_{\frac 12}((0,0))\setminus B_{\frac 1a}((0,0))$, 
$R^1_a=B_{\frac 12}((1,0))\setminus B_{\frac 1a}((1,0))$ and $R^\infty=R_a\setminus(R^0_a\cup R^1_a)$.

By symmetry the integral over $R^0_a$ is the same as that over $R^1_a$ and we have
\[\begin{split}
\iint_{R^0_a} \frac{du \ dv }{(u^{2}+ v^{2})^{2H} \ ((u-1)^{2}+ v^{2})^{2H}}
&\le
 \iint_{R_a^0}
  \frac{du \ dv }{(u^{2}+ v^{2})^{2H} (\frac 12)^{4H}}
 \\&=
 2^{4H} 2\pi \int_{1/a}^{1/2}\frac{r\ dr}{r^{4H}} 
 \\&\le \begin{cases}  \frac{2^{8H-3}\pi}{\frac 12-H}& \text{ if } \frac 14<H<\frac 12\\
                       8\pi \ln a & \text{ if } H=\frac 12\\
                                 \frac{2^{4H-1}\pi}{H-\frac 12}a^{4H-2}& \text{ if } H> \frac 12
 \end{cases}
\end{split}\]
The integral over $R^\infty$ is  a constant $\tilde C_H$ independent of $a$, when $H>\frac 14$. 
Summing up, we get that there exist positive constants $M$ and  $M_H$ such that
\[
I_1(a)\le
\begin{cases}
M_H \frac{1}{a^{8H-2}} &  \text{ if } \frac 14< H<\frac 12\\
M \frac{\ln a}{a^{2}} &  \text{ if } H=\frac 12\\
M_H\frac{1}{a^{4H}} &   \text{ if } H> \frac 12
\end{cases} 
\]
Returning to the notation with $k$ we get our result.
\hfill\qed%\smartqed

%%%%%%%%%%%%LEMMA 2

In the next Lemma we consider a restricted range for $H$; the
assumption on $\rho$ is a restriction of that  in \eqref{cond-sopra}.
\begin{lemma}\label{lemma2}
We are given $\frac 12<H<1$.
Let us assume   $-1<\rho<2(H-1)$.
\\
Then, for any $k \in \zo$ the series
\[
\sum_{h\in\zo, h\neq k} \frac {|h|^{2\rho+2}}{|h|^{4H} |k-h|^{4H}}
\]
converges
and its sum $S_2(k)$ is bounded by
\[
 C_H \frac 1{|k|^{4H-2\rho-2}}
\]
for a suitable constant $C_H$.
\end{lemma}
\proof
By assumption we have $4H+4H-2\rho-2>2$ and therefore the series is
convergent.
To estimate its sum we proceed as in the proof of the previous Lemma.
First the series  can be estimated by an  integral and it is enough to evaluate
\[
I_2(a)=\iint_{\tilde D_a} \frac {dx\ dy}{[x^{2}+ y^{2}]^{2H-\rho-1} \ [(x-a)^{2}+ y^2]^{2H}}
\]
for $a\ge 2$.
By the change of variables $u=\frac xa$ and $v=\frac ya$ we obtain
\[
I_2(a)=
\frac{1}{a^{2(4H-\rho-1)}} \iint_{R_a}
 \frac{a^2}{(u^{2}+ v^{2})^{2H-\rho-1} \ ((u-1)^{2}+ v^{2})^{2H}} \ du \ dv .
\]
We split the integral over $R_a$ into three parts by setting
$R_a=R^0_a\cup R^1_a\cup R^\infty$ with disjoint unions:
\begin{itemize}
\item
over $R^0_a=B_{\frac 12}((0,0))\setminus B_{\frac 1a}((0,0))$
\[\begin{split}
\iint_{R^0_a}
 \frac{du \ dv }{(u^{2}+ v^{2})^{2H-\rho-1} \ ((u-1)^{2}+ v^{2})^{2H}} 
&\le
2^{4H} \iint_{R^0_a}
 \frac{du \ dv }{(u^{2}+ v^{2})^{2H-\rho-1} }
\\&=
2^{4H} 2\pi \int_{1/a}^{1/2}\frac{r\ dr}{r^{2(2H-\rho-1)}}
\end{split}
\]
By assumption we get that 
$2(2H-\rho-1)-1>1$; therefore the latter integral  is bounded by
\[
C_H a^{4H-2\rho-4}
\]
\item
over $R^1_a=B_{\frac 12}((1,0))\setminus B_{\frac 1a}((1,0))$ we
proceed as in the proof of the previous Lemma and get that the
integral over $R^1_a$ is bounded by 
\[
C_H a^{4H-2}
\]
\item
over $R^\infty=R_a\setminus(R^0_a\cup R^1_a)$: this integral is
bounded, uniformly in $a$.
\end{itemize}
Now we compare the exponents of $a$;
since  $\rho>-1$, we have  $4H-2\rho-4 <4H-2$.  Summing the three
contributions and noticing that $4H-2>0$ we conclude
that the integral over the region $R_a$ is bounded by $C_H a^{4H-2}$ for
all $a\ge 2$.
Thus
\[
I_2(a)\le C_H \frac 1{a^{4H-2\rho-2}} \qquad \forall a\ge 2
\]
for a suitable constant $C_H$.
\hfill\qed%\smartqed

%%%%%%%%%%%%%%%LEMMA

\begin{lemma}\label{lemma3}
Let $\frac 14<H<1$ and
\begin{align}
& \rho<4H-3& \text{ if } \frac 14<H<\frac 12 \tag{\ref{cond-sotto}}\\
                        &      \rho< 2(H-1)& \text{ if } \frac 12\le H<1 \tag{\ref{cond-sopra}}
\end{align}
Then
\[\sum_{j\in\zo} |j|^{2\rho+2}  \sum_{h \neq j}   \sum_{l \neq j}
   \frac {|h-l|^{2\rho+2}}{|h|^{4H} |l |^{4H}  |h-j|^{4H} |l-j|^{4H}}<\infty.
\]
\end{lemma}
\proof
First let us prove it when $2\rho+2>0$; this is possible only  when $\frac 12<H<1$. 
In this case we have
\begin{equation}\label{h-l}
|h-l|\le 2 |h| |l| 
\end{equation}
which holds for any $h,l \in \zo$. It comes from
\[
1+\frac{|l|}{|h|} \le 1+|l| \le 2|l|  \qquad \forall |h|, |l|\ge 1;
\]
this is equivalent to 
\[
|h|+|l|\le 2 |h| |l|
\]
Thus, by triangle inequality we obtain \eqref{h-l}.

Hence with a positive power we get
\[
|h-l|^{2\rho+2} \le C |h|^{2\rho+2} |l^{2\rho+2}| 
\]

This implies that we study the series
\[
\sum_{j\in\zo} |j|^{2\rho+2}  \sum_{h \neq j}   \sum_{l \neq j}
   \frac {|h|^{2\rho+2}|l|^{2\rho+2}}{|h|^{4H} |l |^{4H}  |h-j|^{4H} |l-j|^{4H}}
   =
\sum_{j\in\zo} |j|^{2\rho+2}  \left(\sum_{h \neq j}  
   \frac {|h|^{2\rho+2}} {|h|^{4H}   |h-j|^{4H}}\right)^2   .
\]
The innner series is estimated by Lemma \ref{lemma2}; thus
\[
\sum_{j\in\zo} |j|^{2\rho+2}  \left(\sum_{h \neq j}  
   \frac {|h|^{2\rho+2}} {|h|^{4H}   |h-j|^{4H}}\right)^2  
\le C_H
\sum_{j\in\zo}  \frac 1{|j|^{2(4H-3\rho-3)}}
\]
The assumption \eqref{cond-sopra} implies that $2(4H-3\rho-3)>2$ and
therefore this latter series is convergent.

Now let us consider the case  $2\rho+2\le  0$. 
We have $|h-l|^{(2\rho+2)}\le 1$.  
Therefore we are left with
\[
\sum_{j\in\zo} |j|^{2\rho+2}  \sum_{h \neq j}   \sum_{l \neq j}
   \frac 1{|h|^{4H} |l |^{4H}  |h-j|^{4H} |l-j|^{4H}}
   =
   \sum_{j\in\zo} |j|^{2\rho+2} 
   \left( \sum_{h \neq j}     \frac 1{|h|^{4H}   |h-j|^{4H} }\right)^2
\]
We handle this contribution according to Lemma \ref{lemma1}. Indeed it is bounded by
\[
\sum_{j\in\zo}  |j|^{2\rho+2} \times 
\begin{cases}
M_H\frac 1{ |j|^{2(8H-2)}} &  \text{ if } \frac 14<  H<\frac 12\\
M\frac {\ln |j|}{ |j|^{4}} &  \text{ if }  H=\frac 12\\
M_H\frac 1{ |j|^{8H}} &  \text{ if } H> \frac 12
\end{cases}
\]
Thus there is convergence if
\[
\begin{cases}
16H-4-2\rho-2>2&  \text{ if } \frac 14<  H<\frac 12\\
4-2\rho-2>2 &  \text{ if }  H=\frac 12\\
8H-2\rho-2>2 &  \text{ if } H> \frac 12
\end{cases}
\]
These conditions are fulfilled under assumptions \eqref{cond-sotto} and \eqref{cond-sopra}.
\hfill\qed%\smartqed

%%%%%%%%%%%%%
\section{The case $m=2$}\label {m2}
We present the proof of Proposition \ref{stimaBz} for the case $m=2$.

First we shall  use many times  that
$\gamma_{h,k}=\gamma_{k-h,k}$ and $\gamma^2_{h,k}\le |k|^2$ for any $h$.
We have
\[\begin{split}
\int\|B(z,z)\|_{\mathcal H^\rho}^4\ &  \mu^H(dz)
\\
&=\int \left(\sum_{k\in\zo} |k|^{2\rho}|B_k(z,z)|^2\right)^2\ \mu^H(dz)
\\&
=\sum_{k\in\zo}\sum_{j\in\zo}  |k|^{2\rho}|j|^{2\rho}
 \int |B_k(z,z)|^2|B_{j}(z,z)|^2\ \mu^H(dz)
\\&=
\sum_{k\in\zo}\sum_{j\in\zo}  |k|^{2\rho}|j|^{2\rho}
\sum_{h,h^\prime\in\zo; h,h^\prime\neq k}
\sum_{l,l^\prime\in\zo; l,l^\prime\neq j}
\\&\qquad
\int \gamma_{h,k} z_h z_{k-h} \gamma_{h^\prime,k} \overline{z_{h^\prime} z_{k-h^\prime}}
 \gamma_{l,j} z_l z_{j-l} \gamma_{l^\prime,j} \overline{ z_{l^\prime} z_{j-l^\prime}}
  \ \mu^H(dz)
\end{split}
\]
Now we consider the series. There are indeed 6 sums but we have to
consider only the non vanishing integrals: $\mu^H$ is the product
of centered Gaussian measures and all odd powers give zero
contribution in the integral with respect to the measure $\mu^H$.
In particular we analyze 
\[
\int  z_h z_{k-h}  \overline{z_{h^\prime}}\ \overline{ z_{k-h^\prime}}
  z_l z_{j-l}  \overline{ z_{l^\prime}}\ \overline{ z_{j-l^\prime}}
  \ \mu^H(dz)
\]
We get different contributions according to the choice of equal
indices.
Let us list all these possible contributions; we choose $h$
equal to some subsequent index but the case leading to $z_h=\overline{z_{k-h}}$
is not possible since $k \neq 0$. Therefore there are 6 possible
cases.
In the following we do not specify that the sum involves indices
belonging to $\zo$ in order to shorten the notation.
\begin{enumerate}
%1
\item When $h=h^\prime$:
 \begin{itemize}
 \item  for $l=l^\prime$ we get
 \[\begin{split}
 \sum_{k}\sum_{j}  |k|^{2\rho}|j|^{2\rho}&
 \sum_{h \neq k}
 \sum_{l\neq j}
 \int \gamma^2_{h,k} |z_h|^2 |z_{k-h}|^2 
  \gamma^2_{l,j} |z_l|^2 |z_{j-l}|^2 \ \mu^H(dz)
 \\&
 =\left( \sum_{k}  |k|^{2\rho} \sum_{h \neq k}
 \int \gamma^2_{h,k} |z_h|^2 |z_{k-h}|^2   \ \mu^H(dz)\right)^2
 \end{split}
 \]
 We know from \eqref{prima-stima} that the quantity in the parenthesis
 is finite under assumptions \eqref{cond-sotto} and \eqref{cond-sopra}.
 \item  for $l=j-l^\prime$  we get as before.
 \end{itemize}
%2
\item   When $h=k-h^\prime$:
 \begin{itemize}
 \item for $l=l^\prime$  we get as before;
 \item   for $l=j-l^\prime$  we get as before.
\end{itemize}
%3
\item   When $h=-l$   we get
\[
\sum_{k}  \sum_{j}  |k|^{2\rho}|j|^{2\rho}
\sum_{h \neq k}\sum_{h^\prime \neq k} \sum_{l^\prime\neq j}
 \gamma_{h,k} \gamma_{h^\prime,k} \gamma_{-h,j}\gamma_{l^\prime,j}
\int |z_h|^2 z_{k-h}  \overline{z_{h^\prime}}\ \overline{z_{k-h^\prime}}
   z_{j+h}  \overline{ z_{l^\prime}}\ \overline{ z_{j-l^\prime}}
  \ \mu^H(dz)
  \]
  This gives
 \begin{itemize}
 \item for  $k-h=h^\prime$ or $h=h^\prime$
  \[\begin{split}
  \sum_{k} & \sum_{j}  |k|^{2\rho}|j|^{2\rho}
\sum_{h \neq k,-j} \sum_{l^\prime\neq j}
\int \gamma_{h,k}^2 |z_h|^2 |z_{k-h}|^2   \overline{z_{h}}
 \gamma_{-h,j}  z_{j+h} \gamma_{l^\prime,j} \overline{ z_{l^\prime}}\ \overline{ z_{j-l^\prime}}
  \ \mu^H(dz)
\\&
=\sum_{k}\sum_{j}  |k|^{2\rho}|j|^{2\rho}
\sum_{h \neq k,-j}  \gamma_{h,k}^2 (\gamma_{-h,j}^2+\gamma_{-h,j}  \gamma_{h+j,j})
\int  |z_h|^4 |z_{k-h}|^2  |z_{j+h}|^2  \ \mu^H(dz)
\\&\le C
\sum_{k}\sum_{j}  |k|^{2\rho+2}|j|^{2\rho+2}
\sum_{h \neq k,-j} \frac 1{|h|^{8H} |k-h|^{4H}  |j+h|^{4H} }
\\&= C
\sum_{h} \frac 1{|h|^{8H} }
\left(\sum_{k \neq h}\frac {|k|^{2\rho+2}}{ |k-h|^{4H} }\right)^2
\end{split}
\]
The inner series converges
for $\rho<2(H-1)$ (i.e. $2\rho+2-4H<-2$); from \eqref{cond-sotto} and
\eqref{cond-sopra} we know that this condition is fulfilled.
Moreover its sum depends on
$h$ in such a way that it vanishes when $|h| \to \infty$; hence 
it is bounded. Thus we are left with the convergence of 
the first series $\sum_{h} |h|^{-8H} $, which holds for any $H>\frac 14$.
So this  contribution is finite.
\item for $k=-j$
\[\begin{split}
\sum_{k} & |k|^{4\rho}
\sum_{h \neq k} \sum_{h^\prime \neq k} \sum_{l^\prime\neq -k}
\int \gamma_{h,k} |z_h|^2 |z_{k-h}|^2 \gamma_{h^\prime,k} \overline{z_{h^\prime}}\ \overline{z_{k-h^\prime}}
 \gamma_{-h,-k}   \gamma_{l^\prime,-k} \overline{ z_{l^\prime}}\ \overline{ z_{-k-l^\prime}}
  \ \mu^H(dz)
\\&=
\sum_{k}  |k|^{4\rho}
\sum_{h \neq k} \sum_{h^\prime \neq k} 
\gamma_{h,k} \gamma_{h^\prime,k} \gamma_{-h,-k} (\gamma_{-h^\prime,-k} +\gamma_{h^\prime-k,-k})
\int  |z_h|^2 |z_{k-h}|^2  |z_{h^\prime}|^2 |z_{k-h^\prime}|^2 \ \mu^H(dz)
\\&\le C
\sum_{k}  |k|^{4\rho+4}\sum_{h \neq k} \sum_{h^\prime \neq k}
\frac 1{|h|^{4H} |k-h|^{4H}  |h^\prime|^{4H}  |k-h^\prime|^{4H}}
\\&=C
\sum_{k}  |k|^{4\rho+4}\left(\sum_{h \neq k}
\frac 1{|h|^{4H} |k-h|^{4H}}\right)^2
\end{split}
\]
We handle this contribution according to Lemma \ref{lemma1}. Indeed it
is bounded by the series
\[
\sum_{k}  |k|^{4\rho+4} \times 
\begin{cases}
M_H\frac 1{ |k|^{8H-2}} &  \text{ if } \frac 14<  H<\frac 12\\
M\frac {\ln |k|}{ |k|^{2}} &  \text{ if }  H=\frac 12\\
M_H\frac 1{ |k|^{4H}} &  \text{ if } H> \frac 12
\end{cases}
\]
Thus there is convergence if
\[
\begin{cases}
8H-2-4\rho-4>2&  \text{ if } \frac 14<  H<\frac 12\\
2-4\rho-4>2 &  \text{ if }  H=\frac 12\\
4H-4\rho-4>2 &  \text{ if } H> \frac 12
\end{cases}
\]
These conditions are fulfilled under assumptions \eqref{cond-sotto} and \eqref{cond-sopra}.
\item for $k-h=l^\prime$
\[\begin{split}
 \sum_{j}& |j|^{2\rho}   \sum_{h \neq -j}   \sum_{l^\prime \neq -h,j} 
\sum_{h^\prime \neq h+l^\prime}   |h+l^\prime|^{2\rho}
\gamma_{h,h+l^\prime} \gamma_{h^\prime,h+l^\prime} \gamma_{-h,j} \gamma_{l^\prime,j}
\int  |z_h|^2 |z_{l^\prime}|^2  \overline{z_{h^\prime}}\ \overline{z_{h+l^\prime-h^\prime}}
  z_{j+h}  \ \overline{ z_{j-l^\prime}}
  \  \mu^H(dz)
  \\&=
 \sum_{j} |j|^{2\rho}     \sum_{h \neq -j}   \sum_{l^\prime \neq -h,j}
   |h+l^\prime|^{2\rho}
\gamma_{h,h+l^\prime}( \gamma_{h^\prime,h+l^\prime}+\gamma_{l^\prime-j,h+l^\prime}) \gamma_{-h,j} \gamma_{l^\prime,j}
\int  |z_h|^2 |z_{l^\prime}|^2  |z_{j+h}|^2 |z_{l^\prime-j}|^2
\  \mu^H(dz)
 \\&\le C
  \sum_{j} |j|^{2\rho+2}  \sum_{h \neq -j}   \sum_{l^\prime \neq -h,j}
   \frac {|h+l^\prime|^{2\rho+2}}{|h|^{4H} |l^\prime |^{4H}  |j+h|^{4H} |l^\prime-j|^{4H}}
  \end{split}
\]
The convergence of this sums is done in Lemma \ref{lemma3}
in the Appendix.
\item for $k-h=j-l^\prime$, we proceed as in the previous case.
\end{itemize}
%4
\item   when $h=l-j$ we get
\begin{multline*}
\sum_{k}\sum_{j}  |k|^{2\rho}|j|^{2\rho}
\sum_{h\neq k,-j}\sum_{h^\prime\neq k}
\sum_{l^\prime\neq j} 
\gamma_{h,k} \gamma_{h^\prime,k}  \gamma_{h+j,j}   \gamma_{l^\prime,j}
\\
\int  |z_h|^2 z_{k-h}  \overline{z_{h^\prime}}\
 \overline{z_{k-h^\prime}} z_{h+j} \overline{ z_{l^\prime}}\ \overline{z_{j-l^\prime}}
  \ \mu^H(dz)
\end{multline*}
Estimating the $\gamma$'s (that is $|\gamma_{\cdot,k}|\le C|k|$) we get the same computations as in step 3.
%5
\item   when $h=l^\prime$ 
we get the same computations as in step 3. 
%6
\item   when $h=j-l^\prime$  
we get the same computations as in step 3.
\end{enumerate}

This concludes the proof of  \eqref{stimaB(z,z)} for $m=2$.

\end{document}